\documentclass{article}
\usepackage{latexsym}
\begin{document}

\begin{center}
{\bf  Models of {\bf Z}-Orbits of 
            Unitary in Indefinite Inner Product Spaces Operators }
\bigskip

 {\sc Sergej A. Choroszavin } 
\bigskip

\bigskip

 Email: { choroszavin@narod.ru }

\end{center}

\begin{abstract}

 Given a lineal $H_0$ and $x_0\in H_0$ and a linear injective 
 $U_0: H_0 \to H_0$ 
 such that all ${U_0}^N \,,\quad N \in {\bf Z}$ exist 
 and all 
 $\{{U_0}^N x_0 \,|\,N \in {\bf Z} \}$
 are linearly independent, anyone can define 
 on 
 $L_{x_0}:=span\{{U_0}^{\bf Z} x_0\} 
           := span\{{U_0}^N x_0 \,|\,N \in {\bf Z} \}$ 
 a (pre)hilbert scalar product such that $U_0$ becomes a unitary operator.

 The problem under consideration is: Suppose there is specified an 
 indefinite inner product $\{,\}_0$ on $H_0$ and $U_0$ is 
 a $\{,\}_0$-unitary operator. Can one introduce a 
 (pre)hilbert topology on $L_{x_0}$ so that after completion 
 and possible extension
 the resulting $\{,\}_{ext}$ is continuous,
 the resulting $U_{ext}$ is $\{,\}_{ext}$-unitary and there exists a pair 
 $L_{+}, L_{-}$ mutually 
 $\{,\}_{ext}$ orthogonal, maximal strictly positive and respectively 
 negative subspaces, so that 
 they are $U_{ext}$-invariant? More generally, can one construct
 a sequence (chain) of transformations of the type 
$$
  \mbox{ restrict } \to 
  \mbox{ change topology } \to 
  \mbox{ make completion } \to 
  \mbox{ extend } 
  \quad 
$$
$$
  \to 
  \mbox{ restrict } \to 
  \cdots 
$$
 with the same result?
 (And, as a result, after some transformations, which are natural in the field of 
  indefinite inner  product spaces, 
  $U_{ext}$ will become usual Hilbert space unitary operator).

 For a relatively wide class of pairs "operator, inner product" positive 
 solutions proposed.
\end{abstract}

\par\addvspace{\bigskipamount}\par\noindent

\par\addvspace{\bigskipamount}\par\noindent

\newpage

\subsection*{0. Introduction }
%
\par\noindent 
 Given a lineal $H_0$ and $x_0\in H_0$ and a linear injective 
 $U_0: H_0 \to H_0$ 
 such that all ${U_0}^N \,,\quad N \in {\bf Z}$  exist,
 and all 
 $\{{U_0}^N x_0 \,|\,N \in {\bf Z} \}$
 are linearly independent, 
 anyone can define 
 on 
 $L_{x_0}:=span\{{U_0}^{\bf Z} x_0\} 
           := span\{{U_0}^N x_0 \,|\,N \in {\bf Z} \}$ 
 a (pre)hilbert scalar product such that $U_0$ becomes a unitary operator.
 But what about the case where $H_0$ has an additional structure and some 
 properties of $H_0$ are expressed in terms of that structure?
 Namely, what about the case where on $H_0$ is defined an indefinite inner 
 product, say $\{\cdot , \cdot \}_0$, and $U_0$ is $\{,\}_0$-unitary at least 
 on $span\{{U_0}^{\bf Z} x_0\}$ :
$$
    \{U_0 x, U_0 y \}_0 = \{x, y \}_0  \,, \quad 
    \{U_0^{-1} x, U_0^{-1} y \}_0 = \{ x, y \}_0 
$$
for all $x,y \in span\{{U_0}^{\bf Z} x_0\}$, what about such a case?
 If $\{x_0, U_0^N x_0 \}_0$ as a function of $N \in {\bf Z}$ can be presented 
 as linear combination of positive defined functions of $N$, then one can 
 define on $span\{{U_0}^{\bf Z} x_0\}$ a (pre)hilbert scalar product
 (may be degenerated),
 say $(,)$, such that {\bf $U_0$ becomes $(,)$-unitary} and 
 {\bf simultaneously $\{,\}$ becomes continuous} with respect to the 
 corresponding (pre)hilbert norm (may be seminorm).
 If not, if the presentation, mentioned above, is impossible,
 then such $(,)$ is impossible. But what is possible?
\par\addvspace{2\medskipamount}\par\noindent 
 In this paper we will show that at least the following is possible:
\par\addvspace{2\medskipamount}\par\noindent
 Under some natural conditions, one can define a prehilbert scalar product,
 say $(,)$, on $span\{{U_0}^{\bf Z} x_0\}$, so that $\{,\}_0$ turned out to be 
 continuous and after the corresponding 
 completion of  $span\{{U_0}^{\bf Z} x_0\}$ , denote it here by  
 ${L_{x_0}}_\sim$, this space 
 ${L_{x_0}}_\sim$, indefinite inner product $\{,\}$ and $U_0$ 
 can be extended to the Hilbert 
 space $H_{ext}$, continuous indefinite inner product $\{,\}_{ext}$ 
 on $H_{ext}$  and 
 a $\{,\}_{ext}$-unitary operator $U_{ext}: H_{ext}\to H_{ext} $ 
 with the properties: 
\par\addvspace{2\medskipamount}\par\noindent  
 there exist two lineals 
$$
 L_{+} \,,\, L_{-} \subset H_{ext}
$$
such that
\par\addvspace{\medskipamount}\par\noindent 
(a) $U_{ext}^{\pm 1}L_{+} = L_{+} $ 
\par\addvspace{\medskipamount}\par\noindent 
(b) $U_{ext}^{\pm 1}L_{-} = L_{-} $ 
\par\addvspace{\medskipamount}\par\noindent 
(c) $\{x,x\}_{ext} > 0 \quad ( x \in L_{+}\backslash \{0\} )$  
\par\addvspace{\medskipamount}\par\noindent 
(d) $\{x,x\}_{ext} < 0 \quad ( x \in L_{-}\backslash \{0\} )$ 
\par\addvspace{\medskipamount}\par\noindent 
(e)  $\{ L_{+} , L_{-}\}_{ext}=\{0\} $  
\par\addvspace{\medskipamount}\par\noindent 
(f)  $L_{+} \cap  L_{-} = \{0\}$ 
\par\addvspace{\medskipamount}\par\noindent 
(g) if $\{ L_{+} + L_{-}, x \}_{ext}=0$  then  $x = 0$  
\par\addvspace{\medskipamount}\par\noindent 
(h)  $\overline{ L_{+} + L_{-} }= H $ 
\par\addvspace{\medskipamount}\par\noindent 
 Moreover,   
\par\addvspace{\medskipamount}\par\noindent 
(c') $ \{x,x\}_{ext} > 0 \quad ( x \in \overline{L_{+}}\backslash \{0\} ) $  
\par\addvspace{\medskipamount}\par\noindent
(d') $ \{x,x\}_{ext} < 0 \quad ( x \in \overline{L_{-}}\backslash \{0\} ) $ 
\par\addvspace{\medskipamount}\par\noindent 
(e') $ \{ \overline{L_{+}} , \overline{L_{-}}\}_{ext}=\{0\} $ 
\par\addvspace{\medskipamount}\par\noindent 
(f') $ \overline{L_{+}} \cap  \overline{L_{-}} = \{0\} $ 
\par\addvspace{\medskipamount}\par\noindent 
 Indeed, if in addition $U$ and $U^{-1}$ are bounded, then 
\par\addvspace{\medskipamount}\par\noindent 
(a')  ${U_{ext}}^{\pm 1}\overline{L_{+}} = \overline{L_{+}}$ 
\par\addvspace{\medskipamount}\par\noindent
(b')  ${U_{ext}}^{\pm 1}\overline{L_{-}} = \overline{L_{-}}$
\par\addvspace{\medskipamount}\par\noindent
(overline stands for "closure of").
\par\addvspace{\medskipamount}\par\noindent  
\par\addvspace{\medskipamount}\par\noindent  
 Of course, there exists $x_{0 \, ext} \in H_{ext}$ such that 
\begin{eqnarray*}
 \{ x_{0 \, ext} , U_{ext}^N x_{0 \, ext}\}_{ext}
& = &
 \{ x_0 , U_0^N x_0\}_0  \,, \quad   N \in {\bf Z}
\end{eqnarray*}
 and one of methods of this paper consists in to take a Hilbert space $H_1$, 
 indefinite inner product $\{,\}_1$ with desirable properties (a), (b), and so on,
 and then look for $x_1\in H_1$ such that 
\begin{eqnarray*}
 \{ x_1 , U_1^N x_1\}_1
& = &
 \{ x_0 , U_0^N x_0\}_0  \,, \quad   N \in {\bf Z}
\end{eqnarray*}
 If we have success, we can identify "$1$"-objects with "ext"-objects.
 This is why the paper is called "Models of {\bf Z}-orbits ...".
\par\addvspace{\medskipamount}\par\noindent
\par\addvspace{\medskipamount}\par\noindent
  Logically we cosider the chain of transformaitions  
$$
  \mbox{ restrict } \to 
  \mbox{ change topology } \to 
  \mbox{ make completion } \to 
  \mbox{ extend } 
  \quad 
$$
$$
  \quad\quad 
  \to\mbox{ restrict } \to 
  \mbox{ change topology } \to 
  \mbox{ make completion } \to 
  \mbox{ extend }
$$
$$
  \to\mbox{ restrict } \to 
  \cdots 
$$
 The order of the components of the chain may be changed.
\par\addvspace{\medskipamount}\par\noindent
\par\addvspace{\medskipamount}\par\noindent
 One can say that we want to develop a specific version 
 of the technique of modelling operators (model operators) 
 with the aim to construct the models with desired properties,
 not to prove something predesigned about the original abstract operator.
 The latter is secondary: the central object is model, not abstract operator. 
\par\addvspace{\medskipamount}\par\noindent
\par\addvspace{\medskipamount}\par\noindent
 Now about the "some natural conditions". 
\\
 First, we will mostly suppose that the restriction $U_0$ on $L_{x_0}$ 
 has no non-zero so called neutral invariant subspace.
\\
 Second, to avoid problems with the definition of invariant subspace of 
 unbounded operator, we will in some constructions suppose that   
$$
    | \{x_0, U_0^N x_0\}_0 | \leq M\,e^{a\, |N| }
$$
 for some numbers $M, a$ and all integers $N$ .
\par\addvspace{\medskipamount}\par\noindent
\par\addvspace{\medskipamount}\par\noindent
 Of course we suppose that the reader is familiar with the basics 
 of the theory of operators acting in abstract indefinite inner product spaces.
\par\addvspace{\medskipamount}\par\noindent
\par\addvspace{\medskipamount}\par\noindent
 In this paper the modelling operators are  combinations of bilateral 
 weighted shifts. We only need of the basics of the basics of the theory of 
 bilateral 
 weighted shifts.
\par\addvspace{\medskipamount}\par\noindent 
\par\addvspace{\medskipamount}\par\noindent 
We should warn about some peculiarities of the system of notations
(in this paper).
For some reasons we prefer to define sesquilinear form to be 
linear in the {\bf second } argument, not in the first. 
Another, given an element $f$ of a lineal and a number $\lambda$, 
we consider $f\,\lambda$ as a synonimous of $\lambda \, f$. 
\par\addvspace{\medskipamount}\par\noindent 
\par\addvspace{\medskipamount}\par\noindent 
Now the details.

\newpage

\subsection*{1. Notes on Bilateral Weighted Shifts }
%
\par\noindent 
\par\addvspace{\medskipamount}\par\noindent
\par\addvspace{\medskipamount}\par\noindent
Let 
$\{u_n\}_{n\in {\bf Z}}$ 
 be a bilateral number sequence; we will suppose that 
$u_n \not= 0$
 for all 
$n\in {\bf Z}$.
Let
 $H$ be a separable Hilbert space, complex or real, with a fixed basis 
$\{b_n\}_n$ indexed with $n \in {\bf Z}$, (bilateral basis). 

 In this case let 
$U$ 
 denote the {\bf shift } 
\footnote{ 
 the full name is: the {\bf bilateral weighted shift of 
$\{b_n\}_n$,
 to the right}.
} 
 that is generated by the formula 
$$
 U\,:\,b_n\,\mapsto\,\frac{u_{n+1}}{u_n}\,b_{n+1} \qquad . \eqno(*)
$$

 The general facts we need are these:

{\bf Observation 1.1}

 One constructs the 
$U$ as follows:
 
 One starts extending the instruction
$(*)$
 on the linear span of the 
$\{b_n\}_{n\in {\bf Z}}$ 
 so that the resulting operator becomes linear.
 That extension is unique and defines a linear densely defined 
 operator, which is here denoted by
$U_{min}$,
 and which is closable.
 The closure of
$U_{min}$
 is just the 
$U$.

 Now then, this 
$U$
 is closed and at least densely defined and injective;
 it has dense range and the action of
$U^N$, 
$U^{*-N}$, 
$U^{*N}U^N$, 
$U^{-N}U^{*-N}$ 
 (for any integer $N$) is generated by 
$$
 U^N:b_n\mapsto \frac{u_{n+N}}{u_n}b_{n+N}\,;\qquad 
 U^{*-N}:b_n\mapsto \frac{u_n^*}{u_{n+N}^*}b_{n+N}\,;
$$
$$
 U^{*N}U^N:b_n\mapsto {|\frac{u_{n+N}}{u_n}|}^2 b_n\,;\qquad 
 U^{-N}U^{*-N}:b_n\mapsto {|\frac{u_n}{u_{n+N}}|}^2 b_n\,.
$$

 In particular,
$U^N$
 is bounded just when the number sequence 
$\{|u_{n+N}/u_n|\}_n$
 is bounded.
 Moreover,
$$  \| U^N \| = sup\{ \; |u_{n+N}/u_n| \;  | n \in {\bf Z} \}    $$ 
\\ 
$\Box$

\par\addvspace{\medskipamount}\par\noindent
\par\addvspace{\medskipamount}\par\noindent
This is all what we need from the theory of  bilateral weighted shifts.

\newpage

%
%
\par\addvspace{2\medskipamount}\par\noindent 
\subsection*{2. Some properties of $\widehat{U}=U \oplus {U^*}^{-1}$  
where $U$ is a bilateral weighted shift }
First of all recall that $\widehat{U}=U \oplus {U^*}^{-1}$
(acting in $H\oplus H$) is unitary 
with respect to the indefinite inner product $\{,\}$, given by 
\begin{eqnarray*}
 \{ x\oplus y , x_1\oplus y_1   \} 
& := &
 (x, y_1) + (y, x_1) ;
\end{eqnarray*}
 where $(,)$ stands for Hilbert scalar product (in $H$)
\par\addvspace{2\medskipamount}\par\noindent 
{\bf Remark 2.1.}
\\
 The reader maybe prefer to deal with the complex Hilbert space and 
\begin{eqnarray*}
 \{ x\oplus y , x_1\oplus y_1   \}_s 
& := &
 i( (x, y_1) - (y, x_1)  ;
\end{eqnarray*}
Recall 
\begin{eqnarray*}
 \{ (I\oplus -iI)(x\oplus y) , (I\oplus iI) (x_1\oplus y_1)   \} 
& := &
 i( (x, y_1) - (y, x_1)
\\&& =  \{ x\oplus y , x_1\oplus y_1 \}_s ; 
\end{eqnarray*}
 and  
\begin{eqnarray*}
 (I\oplus -iI)(U \oplus {U^*}^{-1}) (I\oplus iI) 
& = &
 (I\oplus iI) (U \oplus {U^*}^{-1}) (I\oplus -iI) 
\\
& = &
 U \oplus {U^*}^{-1} ;
\end{eqnarray*}
 So, in the complex space case both constructions are logically equivalent.
 On the other hand the former is algebraically more simple than latter, 
 and we prefer to deal with $\{,\}$.
\\
$\Box$

\par\addvspace{2\medskipamount}\par\noindent 
{\bf Observation 2.1.}
\begin{eqnarray*}
\{ \widehat{U}^N(b_0\oplus b_0), \widehat{U}^M(b_0\oplus b_0) \}
&=& 
0 , \quad (M \not = N)
\\
\{ \widehat{U}^N(b_0\oplus b_0), \widehat{U}^N(b_0\oplus b_0) \}
&=& 
2(b_0,b_0) = 2 > 0
\\
\{ \widehat{U}^N(b_0\oplus -b_0), \widehat{U}^M(b_0\oplus -b_0) \}
&=&
0 , \quad (M \not = N)
\\
\{ \widehat{U}^N(b_0\oplus -b_0), \widehat{U}^N(b_0\oplus -b_0) \}
&=& 
-2(b_0,b_0) = -2 < 0
\\
\{ \widehat{U}^N(b_0\oplus b_0), \widehat{U}^M(b_0\oplus -b_0) \}
&=&
0 .
\\
\end{eqnarray*}

\par\addvspace{2\medskipamount}\par\noindent 
{\bf Definition 2.1.}
\begin{eqnarray*}
 L_{+}
& := &
 span\{\widehat{U}^N(b_0\oplus b_0)   | N \in {\bf Z}\}
 \equiv span\{(U^N b_0\oplus {U^*}^{-N} b_0)   | N \in {\bf Z}\}
\\
 L_{-}
& := &
 span\{\widehat{U}^N(b_0\oplus -b_0)   | N \in {\bf Z}\}
 \equiv span\{(U^N b_0\oplus -{U^*}^{-N} b_0)   | N \in {\bf Z}\}
\\
\end{eqnarray*}
\par\addvspace{2\medskipamount}\par\noindent 
 {\bf Theorem 2.1.} 
\par\addvspace{\medskipamount}\par\noindent 
\begin{eqnarray*}
(1)
&&
  \widehat{U}^{\pm 1}L_{+} = L_{+}
\\
(2)
&&
  \widehat{U}^{\pm 1}L_{-} = L_{-}
\\
(3)
&&
b_N \oplus 0 \in L_{+} + L_{-}   \quad ( N \in {\bf Z} )
\\
&&
0 \oplus b_N \in L_{+} + L_{-}  \quad ( N \in {\bf Z} )
\\
(4)
&&
  \mbox{if } \{ L_{+} + L_{-}, x \}=0 \mbox{ then } x = 0 
\\
(5)
&&
  \overline{ L_{+} + L_{-} }= \widehat{ H  }
\\
(6)
&&
  \{ L_{+} , L_{-}\}=\{0\} 
\\
(7)
&&
 \{x,x\} > 0 \quad ( x \in L_{+}\backslash \{0\} ) 
\\
(8)
&&
 \{x,x\} < 0 \quad ( x \in L_{-}\backslash \{0\} ) 
\\
(9)
&&
  L_{+} \cap  L_{-} = \{0\}
\\
\end{eqnarray*}
 Moreover,   
\begin{eqnarray*}
(6') 
&&
 \{ \overline{L_{+}} , \overline{L_{-}}\}=\{0\} 
\\
(7')
&&
 \{x,x\} > 0 \quad ( x \in \overline{L_{+}}\backslash \{0\} ) 
\\
(8')
&&
 \{x,x\} < 0 \quad ( x \in \overline{L_{-}}\backslash \{0\} ) 
\\
(9')
&&
 \overline{L_{+}} \cap  \overline{L_{-}} = \{0\}
\\
\end{eqnarray*}
 If in addition $U$ and $U^{-1}$ are bounded, then 
\begin{eqnarray*}
(1')
&&
  \widehat{U}^{\pm 1}\overline{L_{+}} = \overline{L_{+}}
\\
(2')
&&
  \widehat{U}^{\pm 1}\overline{L_{-}} = \overline{L_{-}}
\end{eqnarray*}
\par\addvspace{2\medskipamount}\par\noindent 
{\bf Proof} is straightforward.

\par\addvspace{2\medskipamount}\par\noindent 
{\bf Definition 2.2.}
\begin{eqnarray*}
 L_{g\oplus h}
& = & 
span \{ U^N g \oplus {U^*}^{-N} h  |
\quad N = 0, \pm 1, \pm 2, \pm 3, \ldots 
 \}
\end{eqnarray*}
\\
$\Box$

\par\addvspace{2\medskipamount}\par\noindent 
{\bf Definition 2.3.}
\par\addvspace{\medskipamount}\par\noindent 
Given $f \in H $ we write 
$$
  f(n) := (b_n, f) \,;
$$
 Conversely, given $\{f(n)\}_n \in l_2(\bf Z)$ we put 
$$
     f := \sum_{n\in {\bf Z}} b_n(b_n, f) \in H  \,;
$$
\\
$\Box$

\par\addvspace{2\medskipamount}\par\noindent 
{\bf Theorem 2.2.}
\par\addvspace{\medskipamount}\par\noindent 
If 
$$
\frac{\overline u_{-N}}{\overline u_0} f_{0}(-N)
       = \overline{f_{0}(N)}\frac{u_{N}}{u_0}
$$
then $$ \{ L_{b_0 \oplus f_0} , L_{b_0 \oplus -f_0} \} =  \{0\} $$
\par\addvspace{\medskipamount}\par\noindent 
{\bf Proof.}
\\
\begin{eqnarray*}
\{b_0\oplus -f_{0} , \widehat{ U }^N ( b_0\oplus f_{0} ) \}
& = & \{b_0\oplus -f_{0}, U^N b_0\oplus {U^*}^{-N} f_{0} \}
\\
& = & (b_0,{U^*}^{-N} f_{0} ) - (f_{0} , U^N b_0 )
\\
& = & ({U}^{-N} b_0, f_{0} ) - (f_{0} , U^N b_0 )
\\ 
& = & \frac{\overline u_{-N}}{\overline u_0} f_{0}(-N)
       - \overline{f_{0}(N)}\frac{u_{N}}{u_0}
\\ 
& = &  0 
\quad ,
\\\quad N = 0, \pm 1, \pm 2, \pm 3, \ldots 
\end{eqnarray*}
\\
$\Box$

\par\addvspace{2\medskipamount}\par\noindent 
 {\bf Theorem 2.3.}
\par\addvspace{\medskipamount}\par\noindent 
 Suppose 
\begin{eqnarray*}
 \{ L_{b_0 \oplus f_0} + L_{b_0 \oplus -f_0} , g_1\oplus h_1  \}
& = & 
 \{0\} 
\end{eqnarray*}
Then 

\begin{eqnarray*}
h_1 & = &
 0
\end{eqnarray*}
and 
\begin{eqnarray*}
( f_0 , {U}^{N} g_1 )
& = &
 0
 \,\quad N = 0, \pm 1, \pm 2, \pm 3, \ldots 
\end{eqnarray*}

\par\addvspace{2\medskipamount}\par\noindent 
{\bf Proof}
\begin{eqnarray*}
\{ \widehat{ U }^N ( b_0\oplus \pm f_{0} ) , g_1\oplus h_1  \}
& = &
 0
\\&& \quad N = 0, \pm 1, \pm 2, \pm 3, \ldots 
\end{eqnarray*}
  Then 
\begin{eqnarray*}
\{ \widehat{ U }^N ( b_0\oplus 0 ) , g_1\oplus h_1 \}
& = &
 0
\\&& \quad N = 0, \pm 1, \pm 2, \pm 3, \ldots 
\end{eqnarray*}
  and 
\begin{eqnarray*}
\{ \widehat{ U }^N ( 0 \oplus f_0 ) , g_1\oplus h_1 \}
& = &
 0
\\&& \quad N = 0, \pm 1, \pm 2, \pm 3, \ldots 
\end{eqnarray*} 
 Now obtain 
\begin{eqnarray*}
 0
& = &
\{ \widehat{ U }^N ( b_0\oplus 0 ) , g_1\oplus h_1 \}
\\
& = &
\{  ( { U }^N b_0\oplus 0 ) , g_1\oplus h_1 \}
\\ 
& = &
( { U }^N b_0 ,  h_1 )
\\&& \quad N = 0, \pm 1, \pm 2, \pm 3, \ldots 
\end{eqnarray*}
\begin{eqnarray*}
 0
& = &
\{ \widehat{ U }^N ( b_0\oplus 0 ) , g_1\oplus h_1 \}
=
\{  ( { U }^N b_0\oplus 0 ) , g_1\oplus h_1 \}
=
( { U }^N b_0 ,  h_1 )
\\&& \quad N = 0, \pm 1, \pm 2, \pm 3, \ldots 
\end{eqnarray*}
 Hence $h_1 = 0$.
Now 
\begin{eqnarray*}
 0
& = &
\{ \widehat{ U }^N ( 0 \oplus f_0 ) , g_1\oplus h_1 \}
\\
& = &
\{  ( 0 \oplus {U^*}^{-N}f_0 ) , g_1\oplus h_1 \}
\\
& = &
({U^*}^{-N}f_0  ,  g_1 )
\\
& = &
(f_0  , {U}^{-N} g_1 )
\\&& \quad N = 0, \pm 1, \pm 2, \pm 3, \ldots 
\end{eqnarray*} 
Hence 
\begin{eqnarray*}
 0
& = &
(f_0  , {U}^{N} g_1 )
\\&& \quad N = 0, \pm 1, \pm 2, \pm 3, \ldots 
\end{eqnarray*} 
\\
$\Box$

\par\addvspace{2\medskipamount}\par\noindent 
{\bf Corollary 2.1.}
\par\addvspace{\medskipamount}\par\noindent 
Let $U^{\pm 1}$ be bounded, let $f_0 \not= 0$
and suppose that either 
$\overline{span\{U^N g_1| N \in {\bf Z} \}} = \{0\}$ 
or
$\overline{span\{U^N g_1| N \in {\bf Z} \}} = H$ 
\\
Then 
\begin{eqnarray*}
 \{ L_{b_0 \oplus f_0} + L_{b_0 \oplus -f_0} , g_1\oplus h_1  \}
& = & 
 \{0\} 
\end{eqnarray*}
implies  
\begin{eqnarray*}
g_1\oplus h_1 & = &
0 \oplus 0 
\end{eqnarray*}
and as a result 
\begin{eqnarray*}
 \overline{ L_{b_0 \oplus f_0} + L_{b_0 \oplus -f_0} }
& = & 
 H\oplus H  \,, 
\end{eqnarray*}
\begin{eqnarray*}
  \overline{L_{b_0 \oplus f_0}} \cap  \overline{L_{b_0 \oplus -f_0}}
& = &
 \{0\} \,. 
\end{eqnarray*}
$\Box$

\par\addvspace{2\medskipamount}\par\noindent 
{\bf Remark 2.2}
\par\addvspace{\medskipamount}\par\noindent 
If $U, U^{-1}$ are bounded then 
$$
U\overline{span\{U^N g_1| N \in {\bf Z} \}}
=\overline{span\{U^N g_1| N \in {\bf Z} \}}
$$

\par\addvspace{2\medskipamount}\par\noindent 
{\bf Definition 2.4.}
\par\addvspace{\medskipamount}\par\noindent 
We say that $U$ has the {\bf Property A}, iff 
from 
$f_0\not =0 $
and 
\begin{eqnarray*}
 \{ L_{b_0 \oplus f_0} + L_{b_0 \oplus -f_0} , g_1\oplus h_1  \}
 & = & 
 \{0\} 
\end{eqnarray*}
one can deduce that 
\begin{eqnarray*}
g_1\oplus h_1 & = &
0 \oplus 0 
\end{eqnarray*}
\\
$\Box$

\par\addvspace{2\medskipamount}\par\noindent

\newpage

%
%
\par\addvspace{2\medskipamount}\par\noindent
\subsection*{3. Models of Orbits. Existence and Properties. }
\par\addvspace{2\medskipamount}\par\noindent
\par\noindent 
 {\bf Definition 3.1.}
 Let $H_1, \{, \}_1$ and  $H_0, \{,\}_0$ be indefinite inner product spaces,
 $U_1$ and $U_0$ be respectively  $\{, \}_1$- and $\{, \}_0$-unitary operators
 (bounded or not), and let $x_1 \in H_1$,  $x_0 \in H_0$ 

 We say that the {\bf Z}-orbit 
\begin{eqnarray*}
&&
 \{U_1^N x_1 |  N \in {\bf Z}\}
\end{eqnarray*}
 is a model of the orbit 
\begin{eqnarray*}
&&
 \{U_0^N x_0 |  N \in {\bf Z}\}
\end{eqnarray*}
 iff 
\begin{eqnarray*}
 \{ x_1 , U_1^N x_1\}_1
& = &
 \{ x_0 , U_0^N x_0\}_0  \,, \quad   N \in {\bf Z}
\end{eqnarray*}
\\
In this case we write 
$$
 span \{U_1^N x_1 |  N \in {\bf Z}\} \approx_{\{,\}_1,\{,\}_0}
 span \{U_0^N x_0 |  N \in {\bf Z}\}
$$
 and/or simply 
$$
  \{U_1^N x_1 |  N \in {\bf Z}\} \approx_{\{,\}_1, \{,\}_0}  \{U_0^N x_0 |  N \in {\bf Z}\}
$$
$\Box$
\par\addvspace{2\medskipamount}\par\noindent
 Throughout this section 
\begin{eqnarray*}
 H_1 
& := &
 H\oplus H
\end{eqnarray*}
\begin{eqnarray*}
 \{ x\oplus y , x_1\oplus y_1   \}_1 
& := &
 \{ x\oplus y , x_1\oplus y_1   \} 
\\
& := &
 (x, y_1) + (y, x_1) ;
\end{eqnarray*}
 where Hilbert space $H, (,)$ is the same as in the previous sections   
 and let
\begin{eqnarray*}
 U_1 
& := &
 U\oplus {U^*}^{-1}
\end{eqnarray*}
 for a suitable bilateral shift operator $U$ (bounded or not).  
\par\addvspace{2\medskipamount}\par\noindent
 In addition suppose (mostly thechnical suppositions) that 
$u_n$ is choosen so that 
\begin{eqnarray*}
 \overline u_{-n} 
& = &
 u_n \,, \qquad ( n \in {\bf Z} )
\end{eqnarray*}
 and 
\begin{eqnarray*}
\Bigl\{\frac{1}{u_N}{ \{x_0, U_0^N\ x_0 \}_0 }\Bigr\}_{N \in {\bf Z}}
& \in  &
 l_2({\bf Z})
\end{eqnarray*}

\par\addvspace{2\medskipamount}\par\noindent
 First of all we try 
\begin{eqnarray*}
 x_1 
& := &
 b_0 \oplus f_0
\end{eqnarray*}
 for a suitable $f_0$  
\par\noindent 
\par\noindent 
\par\noindent 
\par\noindent 
\par\noindent 
 So, we consider the next equation, where $f_0$ is to be found.
\begin{eqnarray*}
\{b_0\oplus f_{0}, \widehat{ U }^N ( b_0\oplus f_{0} ) \}
=\{b_0\oplus f_{0}, U^N b_0\oplus {U^*}^{-N} f_{0} \}
& = &
\{x_0, U_0^N\ x_0 \}_0
\\
&&\quad , \quad N = 0, \pm 1, \pm 2, \pm 3, \ldots 
\end{eqnarray*}
 Now we obtain 
\begin{eqnarray*}
(b_0, {U^*}^{-N} f_{0}) + (f_{0}, {U}^{N} b_0) 
& = &
\{x_0, U_0^N\ x_0 \}_0  \quad , \quad N = 0, \pm 1, \pm 2, \pm 3, \ldots 
\end{eqnarray*}
\begin{eqnarray*}
( {U}^{-N} b_0, f_{0}) + (f_{0}, {U}^{N} b_0) 
& = &
\{x_0, U_0^N\ x_0 \}_0  \quad , \quad N = 0, \pm 1, \pm 2, \pm 3, \ldots 
\end{eqnarray*}
\begin{eqnarray*}
(\frac{u_{-N}}{u_0}b_{-N}, f_{0}) + (f_{0},\frac{u_N}{u_0} b_N) 
& = &
\{x_0, U_0^N\ x_0 \}_0  \quad , \quad N = 0, \pm 1, \pm 2, \pm 3, \ldots 
\end{eqnarray*}
\begin{eqnarray*}
\frac{\overline u_{-N}}{\overline  u_0}(b_{-N}, f_{0})
+ \frac{\overline u_N}{\overline u_0} (f_{0}, b_N) 
& = &
\{x_0, U_0^N\ x_0 \}_0  \quad , \quad N = 0, \pm 1, \pm 2, \pm 3, \ldots 
\end{eqnarray*}
\begin{eqnarray*}
\frac{\overline u_{-N}}{\overline  u_0}f_{0}(-N)
+ \frac{\overline u_N}{\overline u_0} \overline{f_{0}(N)} 
& = &
\{x_0, U_0^N\ x_0 \}_0  \quad , \quad N = 0, \pm 1, \pm 2, \pm 3, \ldots 
\end{eqnarray*}
 Note  
\begin{eqnarray*}
\{x_0, U_0^N\ x_0 \}_0
& = &
\overline{\{x_0, U_0^{-N}\ x_0 \}_0}
  \quad , \quad N = 0, \pm 1, \pm 2, \pm 3, \ldots 
\end{eqnarray*}
 Recall,   
\begin{eqnarray*}
  u_N
& = &
\overline{u_{-N}}
  \quad , \quad N = 0, \pm 1, \pm 2, \pm 3, \ldots 
\end{eqnarray*}
  so that 
\begin{eqnarray*}
 ( f_{0}(-N)
+ \overline{f_0(N)}) \frac{u_N}{u_0} 
& = &
\{x_0, U_0^N\ x_0 \}_0  \quad , \quad N = 0, \pm 1, \pm 2, \pm 3, \ldots 
\end{eqnarray*}
  The formal solution is 
\begin{eqnarray*}
  f_{0}(N)
& = &
\frac{1}{2}\frac{\overline u_0}{\overline u_N}\overline { \{x_0, U_0^N\ x_0 \}_0 }
+f_{00}(N)
\quad , \quad N = 0, \pm 1, \pm 2, \pm 3, \ldots 
\end{eqnarray*}
 where 
\begin{eqnarray*}
  \overline f_{00}(-N)
& = &
-f_{00}(N)
\quad , \quad N = 0, \pm 1, \pm 2, \pm 3, \ldots 
\end{eqnarray*}
Fix $f_{00} := 0 $ so that  
\begin{eqnarray*}
  f_{0}(N)
& = &
\overline{f_{0}(-N)}
  \quad , \quad N = 0, \pm 1, \pm 2, \pm 3, \ldots 
\end{eqnarray*}
In this case we have
\begin{eqnarray*}
  \{f_{0}(N)\}_{N \in {\bf Z}}
& \in &
 l_2({\bf Z })
\end{eqnarray*}
 and hence we really have obtained a solution. 
\par\addvspace{2\medskipamount}\par\noindent 
Thus we have obtained 
\par\addvspace{2\medskipamount}\par\noindent 
{\bf Theorem 3.1.}
\par\addvspace{\medskipamount}\par\noindent 
 Choose  
$u_n$ so that 
\begin{eqnarray*}
 \overline u_{-n} 
& = &
 u_n \,, \qquad ( n \in {\bf Z} )
\end{eqnarray*}
 and 
\begin{eqnarray*}
\Bigl\{\frac{1}{u_N}{ \{x_0, U_0^N\ x_0 \}_0 }\Bigr\}_{N \in {\bf Z}}
& \in  &
 l_2({\bf Z})
\end{eqnarray*}
(it is always possible) and put 
\begin{eqnarray*}
  f_{0}(N)
& = &
\frac{1}{2}\frac{\overline u_0}{\overline u_N}\overline { \{x_0, U_0^N\ x_0 \}_0 }
\quad , \quad N = 0, \pm 1, \pm 2, \pm 3, \ldots 
\end{eqnarray*}
 In this case 
\begin{eqnarray*}
  \{f_{0}(N)\}_{N \in {\bf Z}}
& \in &
 l_2({\bf Z })
\end{eqnarray*}
 and if one defines  
$$
     x_1 := b_0 \oplus f_0 
$$
  then    
\begin{eqnarray*}
 \{ x_1 , U_1^N x_1\}_1
& = &
 \{ x_0 , U_0^N x_0\}_0  \,, \quad   N \in {\bf Z}
\end{eqnarray*}
\\
$\Box$
\par\addvspace{4\medskipamount}\par\noindent
  Let us discuss some properties of the obtained solution 
\par\addvspace{4\medskipamount}\par\noindent
Recall the Definition 2.2:
\begin{eqnarray*}
 L_{g\oplus h}
& = & 
span \{ U^N g \oplus {U^*}^{-N} h  |
\quad N = 0, \pm 1, \pm 2, \pm 3, \ldots 
 \}
\end{eqnarray*}
 and recall that by the {\bf Theorem 2.2} we have  
\begin{eqnarray*}
 \{ L_{b_0 \oplus f_0} , L_{b_0 \oplus -f_0} \}
& = & 
 \{0\} 
\end{eqnarray*}
\par\addvspace{2\medskipamount}\par\noindent
 {\bf Theorem 3.2.}
\par\addvspace{\medskipamount}\par\noindent 
 Suppose 
 $H_0$ is a complet normed space with the norm $||\cdot ||_0$
such that $|\{x,y\}_0| \leq ||x||_0||y||_0 \,, \quad (x,y \in H_0)$
and $u_n$ is choosen 
so that 
$$ \sum_{n\in {\bf Z}} |\frac{u_0}{u_n}|^2 ||U_0^{n}\ x_0||_0^2 < \infty $$  
and so that 
$U, U^{-1}$ are bounded 
(it is always possible if $U_0, U_0^{-1}$ are bounded).%
\quad   If in addition 
\begin{eqnarray*}
 \{ L_{b_0 \oplus f_0} + L_{b_0 \oplus -f_0} , g_1\oplus h_1  \}
& = & 
 \{0\} 
\end{eqnarray*}
 then 
$$
 \widetilde {g_1} := \sum_{n\in {\bf Z}} \frac{u_0}{u_n} g_1(n) U_0^{n} x_0
  \in \overline {span\{U_0^{n} x_0 | n\in {\bf Z} \}}
  \subset H_0
$$
  and 
\begin{eqnarray*}
 \{U_0^{N} x_0, \widetilde{g_1} \}_0 
& = &
 0
\\ \quad N = 0, \pm 1, \pm 2, \pm 3, \ldots 
\end{eqnarray*}
\begin{eqnarray*}
 \{U_0^{N} \widetilde{g_1}, \widetilde{g_1} \}_0 
& = &
 0
\\ \quad N  = 0, \pm 1, \pm 2, \pm 3, \ldots 
\end{eqnarray*}
 In particular,
$\overline {span \{U_0^{N} \widetilde{g_1}| n \in {\bf Z}  \} }$
is a $||\cdot||_0$ closed $U_0^{\pm 1}$-invariant $\{,\}_0$-neutral
subspace of $\overline {span \{U_0^{N} x_0 | n \in {\bf Z}  \} } \subset H_0$
(may be $\{0\}$).
\par\noindent 
\par\addvspace{\medskipamount}\par\noindent
 {\bf Proof }
\\
By the Theorem 2.3, 
\begin{eqnarray*}
h_1 & = &
 0
\end{eqnarray*}
\begin{eqnarray*}
( f_0 , {U}^{N} g_1 )
& = &
 0
\\ \quad N = 0, \pm 1, \pm 2, \pm 3, \ldots 
\end{eqnarray*}
 On the other hand 
\begin{eqnarray*}
 g_1
& = &
\sum_{n\in {\bf Z}} g_1(n) b_n 
\end{eqnarray*}
 Hence 
\begin{eqnarray*}
 U^{-N}g_1 
& = &
 \sum_{n\in {\bf Z}} \frac{u_{n-N}}{u_n} g_1(n)b_{n-N} 
\\
&&
 \quad N = 0, \pm 1, \pm 2, \pm 3, \ldots 
\end{eqnarray*}
\begin{eqnarray*}
( f_0 , \sum_{n\in {\bf Z}} \frac{u_{n-N}}{u_n} g_1(n)b_{n-N} ) 
& = &
 0
\\ \quad N = 0, \pm 1, \pm 2, \pm 3, \ldots 
\end{eqnarray*}
\begin{eqnarray*}
  \sum_{n\in {\bf Z}} \frac{u_{n-N}}{u_n}  g_1(n) \overline{f_0(n-N)} 
& = &
 0
\\ \quad N = 0, \pm 1, \pm 2, \pm 3, \ldots 
\end{eqnarray*}
Recall, 
\begin{eqnarray*}
  f_{0}(N) = (b_N , f_0)
& = &
\frac{1}{2}\frac{\overline u_0}{\overline u_N}\overline { \{x_0, U_0^N\ x_0 \}_0 }
\quad , \quad N = 0, \pm 1, \pm 2, \pm 3, \ldots 
\end{eqnarray*}
Hence 
\begin{eqnarray*}
  f_{0}(n-N) = (b_{n-N} , f_0)
& = &
\frac{1}{2}\frac{\overline u_0}{\overline u_{n-N}}
\overline { \{x_0, U_0^{n-N}\ x_0 \}_0 }
\quad , \quad N = 0, \pm 1, \pm 2, \pm 3, \ldots 
\end{eqnarray*}
  Hence 
\begin{eqnarray*}
\sum_{n\in {\bf Z}} 
\frac{u_{n-N}}{u_n} g_1(n)
\frac{1}{2}\frac{ u_0 }{ u_{n-N} } \{x_0, U_0^{n-N}\ x_0 \}_0
& = &
 0
\\ \quad N = 0, \pm 1, \pm 2, \pm 3, \ldots 
\end{eqnarray*}
\begin{eqnarray*}
\sum_{n\in {\bf Z}} 
\frac{u_0}{u_n} g_1(n)
\frac{1}{2}  \{x_0, U_0^{n-N}\ x_0 \}_0 
& = &
 0
\\ \quad N = 0, \pm 1, \pm 2, \pm 3, \ldots 
\end{eqnarray*}
\begin{eqnarray*}
\sum_{n\in {\bf Z}} 
\frac{u_0}{u_n} g_1(n)
\frac{1}{2}  \{U_0^{N} x_0,  U_0^{n}\ x_0 \}_0
& = &
 0
\\ \quad N = 0, \pm 1, \pm 2, \pm 3, \ldots 
\end{eqnarray*}
\begin{eqnarray*}
 \{U_0^{N} x_0, \widetilde{g_1} \}_0 
& = &
 0
\\ \quad N = 0, \pm 1, \pm 2, \pm 3, \ldots 
\end{eqnarray*}
\begin{eqnarray*}
 \{U_0^{N}U_0^{n} x_0, \widetilde{g_1} \}_0 
& = &
 0
\\ \quad N,n  = 0, \pm 1, \pm 2, \pm 3, \ldots 
\end{eqnarray*}
\begin{eqnarray*}
 \{U_0^{N} \widetilde{g_1}, \widetilde{g_1} \}_0 
& = &
 0
\\ \quad N  = 0, \pm 1, \pm 2, \pm 3, \ldots 
\end{eqnarray*}
  Quod erat demonstrandum.
\\
$\Box$
\par\addvspace{2\medskipamount}\par\noindent

\par\addvspace{2\medskipamount}\par\noindent
 {\bf Theorem 3.3.}
\par\addvspace{\medskipamount}\par\noindent 
 Suppose 
 $H_0$ is a complet normed space with the norm $||\cdot ||_0$
such that $|\{x,y\}_0| \leq ||x||_0||y||_0 \,, \quad (x,y \in H_0)$
and $u_n$ is choosen 
so that 
$$P(\lambda)=\sum_{n\in {\bf Z}} \frac{u_0}{u_n} \lambda^{n}$$
as function of complex $\lambda$ is analytic on the spectrum of 
$U_0$
(it is always possible if $U_0, U_0^{-1}$ are bounded)
\par\addvspace{\medskipamount}\par\noindent 
Suppose  
 $U_0 |\overline {span \{U_0^{N} \widetilde{g_1}| n \in {\bf Z}  \} }$ 
 (the restriction of $U_0$ onto 
        $\overline {span \{U_0^{N} \widetilde{g_1}| n \in {\bf Z}  \} }$)
has no non-zero neutral invariant subspace.

\par\addvspace{\medskipamount}\par\noindent 
Then  
$span\{U_0^N x_0| N \in{\bf Z}\}$ 
is finite dimensional 
or 
$U$ has the Property A.
\par\addvspace{2\medskipamount}\par\noindent
 {\bf Proof.}
\par\addvspace{\medskipamount}\par\noindent 
Let 
\begin{eqnarray*}
 \{ L_{b_0 \oplus f_0} + L_{b_0 \oplus -f_0} , g_1\oplus h_1  \}
& = & 
 \{0\}  \quad\quad (*)
\end{eqnarray*}
 Then 
$$
   \{g_1(n)\}_{n\in {bf Z}} \in l_2({\bf Z})
$$
 Then 
$$P_1(\lambda)=\sum_{n\in {\bf Z}} g_1(n)\frac{u_0}{u_n} \lambda^{n}$$
as function of complex $\lambda$ is analytic on the spectrum of 
$U_0$
\\
Then 
$$P_1(U)=\sum_{n\in {\bf Z}} g_1(n)\frac{u_0}{u_n} U_0^{n}$$
 exists.
\\
Then 
$$P_1(U)x_0=\sum_{n\in {\bf Z}} g_1(n)\frac{u_0}{u_n} U_0^{n}x_0$$
 exists.
\\
Then
$$ \sum_{n\in {\bf Z}} |\frac{u_0}{u_n}|^2 ||U_0^{n}\ x_0||_0^2 < \infty $$  
Then
$$P_1(U)x_0=\widetilde {g_1}$$
\\
Supposition about neutral invariant subspaces implies $\widetilde {g_1}=0$
Thus $(*)$ implies 
$$ P_1(U)x_0=0 $$
Suppose $(*)$ implies that $P_1(\lambda)$ is identically zero for any $g_1$. 
In this case $(*)$ implies that $g_1=0$ for any $g_1$ 
and hence $U$ has the Property A.
Now suppose that there exists $g_1$ such that $P_1(\lambda)$ is not 
identically zero. In this case, if for a $\lambda_0$ we have 
$P_1(\lambda_0)=0$, then there exists a natural $n_0$ such that 
$(U_0 -\lambda_0 I)^{n_0}x_0 =0$. Hence 
$span\{U_0^N x_0| N \in{\bf Z}\}$ 
is finite dimensional. 
Finally, if $P_1(\lambda)\not=0$ nowhere, then $P_1(U_0)$ is injective,
hence $x_0=0$ and $span\{U_0^N x_0| N \in{\bf Z}\}$ 
is again finite dimensional. 
\\
$\Box$
\par\addvspace{2\medskipamount}\par\noindent
Thus, we have constructed the following chain of restrictions-extensions:   
$$
 H_0 \supset L_{x_0} \approx_{\{,\}_0, \{,\}} L_{b_0 \oplus f_0} \subset 
 L_{b_0 \oplus f_0} \dot{+}_{\{,\}} L_{b_0 \oplus -f_0}
 \subset_{dense} \overline{ L_{b_0 \oplus f_0} \dot{+}_{\{,\}} L_{b_0 \oplus -f_0}}
$$
$$
= \overline{L_{+} \dot{+}_{\{,\}} L_{-}}  \qquad (chain \; A )
$$
Here and in the following $A\dot{+}B$ stands for $A+B$ with indicating 
$A\cap B =\{0\}$ 
and  $A{+}_{\{,\}}B$ stands for $A+B$ with indicating $\{A,B\} =\{0\}$. 
\par\addvspace{2\medskipamount}\par\noindent
\subsection*{4. Operator $U_0 \oplus U_0$ }
\par\addvspace{2\medskipamount}\par\noindent
In this subsection we suppose that all ${U_0}^N x_0 \quad (N \in {\bf Z})$ 
are linearly independent and the sequence $\{u_n\}_n$ is choosen so that 
$U$ has the Property A. 
\par\addvspace{2\medskipamount}\par\noindent
{\bf Observation 4.1.}
\begin{eqnarray*}
\{b_0\oplus -f_{0}, \widehat{ U }^N ( b_0\oplus -f_{0} ) \}
&=&\{b_0\oplus -f_{0}, U^N b_0\oplus -{U^*}^{-N} f_{0} \}
\\
& = &
-\{x_0, U_0^N\ x_0 \}_0  \quad , 
\\
&&\quad N = 0, \pm 1, \pm 2, \pm 3, \ldots 
\end{eqnarray*}
$\Box$
\\
This motivates the following considerations.

\par\addvspace{2\medskipamount}\par\noindent    
{\bf Definition 4.1.}
\\
Put 
\begin{eqnarray*}
 \widetilde H_0
& := &
 H_0\oplus H_0
\\
\|x \oplus y \|_{\sim}^2  
& := &
\|x\|_{0}^2 + \|y \|_{0}^2 \, \quad ( x,y \in H_0 )
\\
 \widetilde U_0
& := &
 U_0\oplus U_0
\\
 \{ x\oplus y , x_1\oplus y_1     \}_{-}
& := &
 \{ x , x_1 \}_0 - \{ y , y_1 \}_0
\\
 L_{x_0}
& := &
 span\{ U^n x_0 | n \in {\bf Z} \}
\\
\overline{ L_{x_0} } &:=& \mbox{ closure of } L_{x_0} 
                      \mbox{ with respect to } \|\cdot\|_0
\\
\widetilde L_{x_0}
& := &
L_{x_0}\oplus L_{x_0}
\\
\end{eqnarray*}
and define a linear map 
\begin{eqnarray*}
 \Omega 
& : &
 L_{x_0} \oplus L_{x_0} \subset \widetilde H_0
  \to 
 L_{b_0\oplus f_0} + L_{b_0\oplus -f_0} \subset \widehat H
\\
\end{eqnarray*}
 as follows:
Firstly put 
\begin{eqnarray*}
 \Omega ( \, U_0^N x_0 \oplus U_0^M x_0 \,)
& :=  &
 \widehat U^N (b_0\oplus f_0) + \widehat U^M  (b_0\oplus -f_0 ) 
\\
&&
   N,M \in {\bf Z} 
\end{eqnarray*}
Then extend this map onto the span of $U_0^N x_0 \oplus U_0^M x_0$ \quad 
$(N,M \in {\bf Z})$ as a linear map.
\par\addvspace{2\medskipamount}\par\noindent
{\bf Observation 4.2.}
\begin{eqnarray*}
 \{
 \Omega ( \, x\oplus y,) ,
 \Omega ( \, x_1\oplus y_1 \,)
 \}
 & =  &
 \{
   x\oplus y, ,
   x_1\oplus y_1 
 \}_{-}
\\
&&
   x, y, x_1, y_1 \in L_{x_0} 
\end{eqnarray*}
\\
Thus, if we take into account only indefinite inner products, we may say that 
$L_{x_0}\oplus L_{x_0}$ is embedded in  
$L_{b_0\oplus f_0} + L_{b_0\oplus -f_0}$, 
we will write it symbolically as 
$$
L_{x_0}\oplus L_{x_0} \subset_{\Omega, \{,\}_0 , \{,\} }
L_{b_0\oplus f_0} + L_{b_0\oplus -f_0}
$$
Recall, $U$ has the Property. As a consequece, 
\\
$$ \overline{ L_{b_0 \oplus f_0} + L_{b_0 \oplus -f_0} } = \widehat H \, \quad (a) $$
and 
$$  L_{b_0 \oplus f_0} \cap L_{b_0 \oplus -f_0} = \{0\}\,; \quad (b)$$
As a result,  
\begin{eqnarray*}
 \Omega 
& : &
L_{x_0} \oplus L_{x_0} \subset \widetilde H_0
  \to 
 L_{b_0\oplus f_0} + L_{b_0\oplus -f_0} \subset \widehat H
\\
\end{eqnarray*}
is  injective (note, surjective by the definition). Hence, 
\begin{eqnarray*}
        \Omega^{-1} 
& : &
 L_{b_0\oplus f_0} + L_{b_0\oplus -f_0} \subset \widehat H
  \to 
 L_{x_0} \oplus L_{x_0} \subset \widetilde H_0
\\
\end{eqnarray*}
exists, and we may identify 
$L_{b_0\oplus f_0} + L_{b_0\oplus -f_0}$                          
and                                                               
$L_{x_0} \oplus L_{x_0}$
as indefinite inner product lineals, write it symbolically:
$$
L_{x_0}\oplus L_{x_0} \approx_{\Omega, \{,\}_0 , \{,\} }
L_{b_0\oplus f_0} + L_{b_0\oplus -f_0}
$$
Thus, if we take into account the relation $(b)$, we may write finally:
$$
L_{x_0}\oplus L_{x_0} \approx_{\Omega, \{,\}_0 , \{,\} }
L_{b_0\oplus f_0} + L_{b_0\oplus -f_0} 
\qquad\qquad\qquad\qquad 
$$
$$
\qquad\qquad 
\subset 
\overline{L_{b_0\oplus f_0} + L_{b_0\oplus -f_0}}
=\widehat{H}
=\overline{L_{+} \dot{+}_{\{,\}} L_{-}} 
\qquad (chain \; B) ;
$$
                                                                                   
\par\noindent 
$\Box$
\par\addvspace{2\medskipamount}\par\noindent    
As a result we have obtained 
\par\addvspace{2\medskipamount}\par\noindent    
{\bf Theorem 4.1.}
\par\addvspace{\medskipamount}\par\noindent 
Suppose $H_0$ is a Hilbert space.
Suppose 
$U_0$ is bounded, has no nonzero $\{,\}_0$ invariant subspace, and $x_0$ is 
a {\bf Z}-cyclic vector of $U_0$, i.e., 
$$
H_0=\overline{ L_{x_0} } 
$$
Then 
\par\addvspace{2\medskipamount}\par\noindent    
$(a)$ 
there exists a (pre)hilbert norm on $L_{x_0}\oplus L_{x_0}$ such that 
$\{,\}_{-}$ and $U_0\oplus U_0$ are continuous and 
can be extended to continuous $\{,\}_{ext}$ and $U_{ext}$ 
onto $H_{ext} = \mbox{ the completion of } L_{x_0}\oplus L_{x_0}$ ;
\\

\par\addvspace{\medskipamount}\par\noindent    
$(b)$ there are two $U_{ext}^{\pm 1}$-invarint (closed) subspaces 
${\cal L}_{+},{\cal L}_{-} \subset \widetilde{ H_0 }$ such that:
\par\addvspace{\medskipamount}\par\noindent    
(1) ${\cal L}_{+}$ is strictly positive;
\par\addvspace{2\medskipamount}\par\noindent    
(2) ${\cal L}_{-}$ is strictly negative;
\par\addvspace{2\medskipamount}\par\noindent    
(3) ${\cal L}_{+}\cap {\cal L}_{-} = \{0\}$ ,
\par\addvspace{2\medskipamount}\par\noindent    
(4) ${\cal L}_{+}+{\cal L}_{-} $ is dense in $\widetilde{H_0}$,
\par\addvspace{2\medskipamount}\par\noindent    
(5) $\{{\cal L}_{+},{\cal L}_{-}\}_{-}=\{0\} $.
\par\addvspace{2\medskipamount}\par\noindent    
(6) as semidefined subspaces, ${\cal L}_{+}$, ${\cal L}_{-}$ both are maximal.
\par\addvspace{2\medskipamount}\par\noindent    
$\Box$
\par\addvspace{2\medskipamount}\par\noindent    
Finally,
\par\addvspace{\medskipamount}\par\noindent    
{\bf Remark 4.1.}
\\
Note, that, in contrast with $(chain\; A)$ in the previous subsection, 
in $(chain\; B)$ all embeddings are dense.
$\Box$
\par\addvspace{2\medskipamount}\par\noindent



\newpage
\bibliographystyle{unsrt}

\end{document}